\documentclass{article} 

\newenvironment{eproof}{\vspace{-2mm}\par 
\par\bgroup\small{\bf Proof.}}{\qed\egroup\par
}
\newcommand{\both}[1]{{\ifmmode{#1}\else{$#1$}\fi}}
\newcommand{\myBox}{\rule{4pt}{10pt}}
\newcommand{\calA}{\both{{\cal A}}}

\newcommand{\calC}{\both{{\cal C}}}

\newcommand{\calH}{\both{{\cal H}}}

\newcommand{\calM}{\both{{\cal M}}}

\newcommand{\calX}{\both{{\cal X}}}
\newcommand{\child}{\both{{\rm Child}}}

\newcommand{\integer}{\both{\bf Z}} 

\newcommand{\pair}[1]{\both{{\langle#1\rangle}}}

\newcommand{\proves}{\vdash}
\newcommand{\qed}{\hfill \myBox}
\newcommand{\ra}{\both{\rightarrow}}

\newcommand{\setii}[2]{\both{{ \left\{\; {\strut#1} \;\left| \; {\strut#2} \right.\;\right\}}}}
\newcommand{\seti}[1]{\both{{\left\{\; #1 \;\right\}}}}

\newcommand{\Frac}[2]{\displaystyle{\frac{\strut #1}{\strut #2}}}

\usepackage[xdvi]{graphics}

\newtheorem{theorem}{Theorem}[section]
\newtheorem{thm}{Theorem}[section]
\newtheorem{prop}[theorem]{Proposition}
 
\newtheorem{lem}[theorem]{Lemma}

\newcommand{\myheadi}[1]{\vspace{1mm}\par\noindent\textbf{#1}}

\newcommand{\comp}{{\mathbf C}}

\newcommand{\myframe}{\both{\mathbf{frame}}}
 
\newcommand{\PD}{{\cal L}}
 
\newcommand{\psdmtree}{{\cal T}}
\newcommand{\psdmposet}{{\cal M}} 
\newcommand{\hypergraph}{\both{{\cal
      HG}raph}}

\newcommand{\olam}{\overline{\lambda}}
\newcommand{\osigma}{\overline{\sigma}}
\title{Higher Dimensional Hypercategories
} 

\author{ Akira Higuchi \\ Hokkaido University\\ 
  {\small \texttt{a-higuti@math.sci.hokudai.ac.jp}} 
\and Hiroyuki Miyoshi \\
Kyoto Sangyo University \\
{\small \texttt{hxm@cc.kyoto-su.ac.jp}}
  \and Toru Tsujishita\\ Hokkaido University\\ 
{\small \texttt{tujisita@math.sci.hokudai.ac.jp} }}
\date{1999.8.3}
\begin{document}
\maketitle
{\footnotesize
\tableofcontents
}\newpage 
\begin{center}
  \textbf{\large Abstract}
\end{center}
{\small
\emph{We introduce higher dimensional hypergraphs, which is a
generalization of Baez-Dolans's opetopic sets and Hermida-Makkai-Power's
mutigraphs. This is based on a simple combinatorial structure
called \emph{shells} and the \emph{formal composites} of pasting diagrams
based on the \emph{closure} of open shells. 
We give two types of graphical representation of higher dimensional
cells which show effectively the relationship of cells of different dimensions. 
Using the hypergraphs, we define strict hypercategories and illustrate its
use by taking Lafont's interaction combinator\cite{lafont} as an example. 
We also give a definition of weak $\omega$-hypercategories and show
that usual category is identified with a special kind of
weak hypercategory as a sample of arguments in our framework.  }

\section*{Introduction}
The purpose of this paper is to give an elementary reformulation of
Baez-Dolan's definition \cite{baez-dolan} of weak higher dimensional
categories from the following points of views:
\textbf{(H1)} In view of the importance of the higher dimensional category theory
  in wide area of mathematical science, it should be formulated elementary
  without using advanced category theory.
\textbf{(H2)} The shapes of cells and pasting diagrams, prepared by the
framework, should be as general as possible, in order not to restrict 
the potential application.
\textbf{(H3)} Graphical representations of cells should be devised which 
describe succinctly the interrelationship among cells of various dimensions.

The point (H1) may not seem realistic in view of the necessity of the
usage of advanced category theory in various formulations
\cite{batanin,tamsamani} but its feasibility is seen clearly in the
Baez-Dolan's formulation which radically reconstructs category theory 
in the following points:
(i) The composites of a pasting diagram of $n$-cells, if they existed at
all, should be explicitly described $n+1$-cells which they call \emph{universal}. 
(ii) The result of composition of $n$-cells depends on the universal $n+1$-cells
, but the results are as it were \emph{isomorphic}. 
(iii) However, ``isomorphism'' is a primitive notion, not defined one as in 
the usual category theory. In fact, in their formulation, the composition and 
isomorphisms are defined mutually recursively\footnote{Note that the usual 
definition of \emph{inverse} makes no sense since the composition is not unique by 
(ii) and the very notion of \emph{identity} does not exist by (iv).}. 
(iv) The existence of cells which play the role of identities are derived from the axiom
on universal cells. See \cite{makkai} for another radical viewpoint of category theory.

The point {(H2)} has importance in view of the fact that the 
``computad'', a sort of 2-graph, whose 2-cells have general polyhedron
shapes (\cite{power,street-computad}) is indispensable in 2-category
theory. However usually they are considered merely as a useful
informal graphical apparatus, which is theoretically unnecessary. For
example, the pasting diagram shown in the upper Figure
\ref{fig:computad} with 2-cells of general form is considered as an
informal representation of the sequences of globular 2-cells drawn in
the lower, which is not unique but equivalent to any other such
sequence by virtue of the interchange law. The upper diagram
seems to be usually more helpful in actual reasoning compared with the 
lower one.

Most of the current proposals of the formulation of higher dimensional
category restrict the shape of cells, e.g. to globular forms
\cite{batanin} or to tree forms
\cite{baez-dolan,power-hermida,miyoshi}. The definition of ``weak
$n$-category'' seems to be considerably simplified by removing such
restriction. For example, we will see that by using the cells of
general form, we can in a sense unify the notion of
\emph{universality} and \emph{balancedness} in the Baez-Dolan's
formulation.

\begin{figure}[htbp]
  \begin{center}
    \leavevmode
    \includegraphics{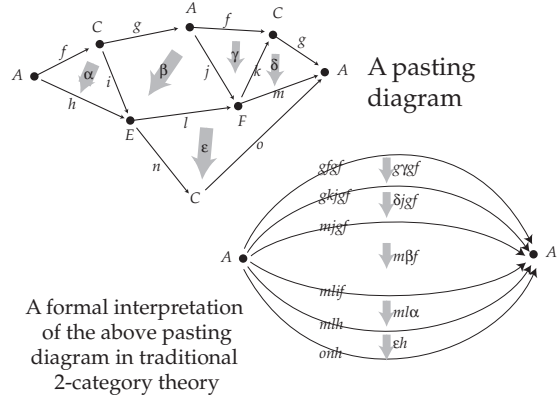}
    \caption{A pasting diagram and one of its formal representations}
    \label{fig:computad}
  \end{center}
\end{figure}
 
Another obstacle of higher dimensional category theory seems to be the
lack of methods of clearly drawing higher dimensional cells which can
conceal unnecessary lower dimensional details appropriately. Note that the direct
generalization of usual 2-pasting diagrams is not feasible for dimensions
higher than 2 and for higher dimensional multicategories.

\myheadi{Outline of the definition} We proceed as follows.
\textbf{(i)} We introduction of a simple combinatorial structure called \emph{shell},
which describes universal shapes of cells and pasting diagrams.
(ii) We define the notion of \emph{labeling} of shells, which substantiate shells, and
certain classes of intermediate substantiation are called \emph{cells, frames} 
and \emph{pasting diagrams}.
We define the \emph{closure} operation and \emph{formal composite} of
pasting diagrams, which is the most important combinatorial apparatus
of this paper.
(iii) A \emph{hypergraph} is defined inductively by substantiating some of the frames of
previous dimension. 
(iv) We show that there is a monad on the category of $n$-hypergraphs.
Strict hypercategories is defined to be the algebra over this monad or
its submonads. 
(v) To show the usefulness of our framework, we describe the process of
graph rewritings in the Lafont's interaction network
as certain 2-dimensional cells in a free strict 2-hypercategory. 
This framework seems to possess a nice connection with graphical language.  
This give also semantics for a programming language \cite{hyper}.
(vi) A hypercategory is an $\omega$-hypergraph endowed with special cells
called universal such that {(iii-a)} pasting diagram of dimension $n$ can be completed to an $n+1$-cell, giving composition like operations which however may not have unique result, {(iii-b)} universal cells are closed under the ``composition'', which play the role both of composition and of equivalence, {(iii-c)}universal cells has inverses with conjugate frames.
(vii) As an example of arguments of our formulation, we prove that a
1-hypercategory of certain type is nothing but a usual category. The
points of proof are the construction of the identity maps and the
proof of the associativity.

\section{Preliminaries}
\paragraph{Trees}
\emph{A tree} is a directed graph with a node called \emph{the root} 
to which there is a unique path from every node other than the root. 
Nodes different from the root are called \emph{general}. 
When there is an edge from $c$ to $p$, we call $p$ \emph{the parent
of} $c$ and $c$ \emph{a child of} $p$. When there is a path from 
a node $x$ to a node $y$, $x$ is called \emph{a descendant of} $y$ 
and $y$ \emph{an antecedent of} $x$. 

Let $T$ be a tree.  The set of the children of a node $p$ is denoted
by $\child_T(p)$ or often simply by $\child(p)$.  A node $x$ and its descendants
form a tree denoted by $T^x$.  

The length of the unique path from a node $x$ to the root is called \emph{the depth of}
$x$. We denote by $T[k]$ the set of all the nodes of depth $k$, so that 
$ T[1]=\child(o_T)$ and $ T[i+1]=\coprod_{j\in T[i]}\child(j)$.
Note that if the node $x$ is of depth $k$ then $T^x[i]\subseteq T[k+i]$.
The sets of nodes of depth $\geq i$ is denoted by $T[[i]]$, namely,
$$ T[[i]]:=\coprod_{j\geq i}T[j].$$  

A tree is called \emph{of height $\leq n$} if the depth of its nodes
are $\leq n$ and \emph{of height $n$} if it is of height $\leq n$ and
there is at least one node of depth $n$.

\paragraph{${\cal C}$-words}
Let $\calC$ be a category. 
\emph{A $\calC$-word} is a family of
objects indexed by a finite set. A $\calC$-word is denoted as 
$w=(w_i)_{i\in |w|}$. An isomorphism $f:w\ra u$ consists of a bijection $|f|:|w|\ra|u|$ and a family of iso's $ f_i:w_i\ra u_{|f|(i)}$ for $i\in |w|$.
A $\calC$-word can be described in various ways. 

For example, a $\calC$-word $w$
with $|w|=\seti{2,3,5,a,b}$ and $w_2=A$, $w_3=B$, $w_5=A$, $w_a=C$, $w_b=A$
can be drawn as follows.
\begin{center}
  \includegraphics{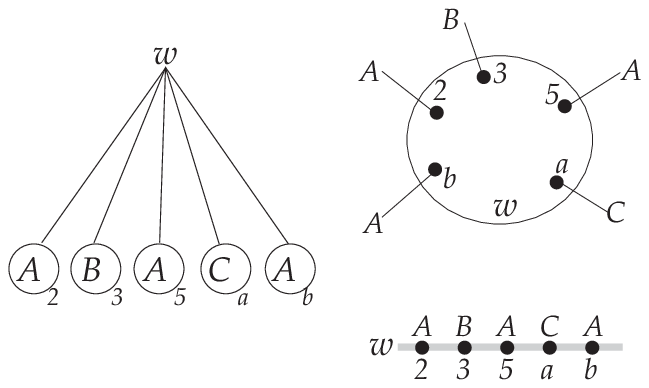}
\end{center}
 
\paragraph{${\cal C}$-trees}
Let $\calC$ be a category. 
\emph{A $\calC$-tree}, or \emph{a tree of $\calC$-objects}, 
is a family of $\calC$-objects indexed by
the nodes of a tree. We write a $\calC$-tree as $w=(w^p)_{p\in T_w}$, where $T_w$ is 
\emph{the underlying tree of} $w$. 
An isomorphism $f:w\ra u$ consists of a tree isomorphism $|f|:T_w\ra T_u$ and isomorphisms 
$$ f_p:w^p\ra u^{|f|p}$$
A $\calC$-tree $w$ defines a $\calC$-word $w[n]:=(w^p)_{p\in T_w[n]}$ for each
nonnegative integer $n$, called \emph{the depth $n$ layer of } $w$. 

\paragraph{${\cal C}$-links}
A finite groupoid is called \emph{a link type} if it has no 
nontrivial endoarrows and its nontrivial isomorphism called \emph{an involution}, 
is composable only with its own inverse and the identities. More explicitly,
a finite groupoid $L$ is a link type if 
\begin{itemize}
\item $L(a,a)=\seti{ id}$,
\item if $a\neq b$, then $|L(a,b)|\leq 1$,
\item if $L(a,b)\bigcap L(b,c)\neq \emptyset$ with $a \neq b\neq c$, then
  $a=c$. 
\end{itemize}

The following is a typical example of a link type,
\begin{center}
  \includegraphics{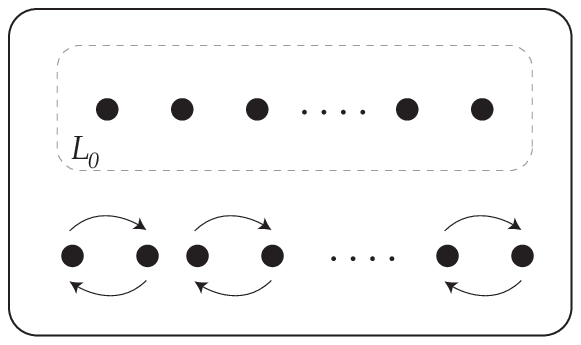},
\end{center}
where the identity arrows are omitted. 

An object of a link type is called \emph{internal} if it is the domain of an involution
and \emph{external} otherwise. External objects form a discrete full subcategory
denoted by $L_{ext}$, which is just a finite set considered as a discrete category.
A link type without external objects is called \emph{closed}. A link type which is
not closed is called \emph{open}.

\emph{A $\calC$-link} $\varphi$ is a functor from a link type
$|\varphi|$ to $\calC$, which we describe explicitly as
$$\varphi=((\varphi_i)_{i\in |\varphi|},(\varphi_m:\varphi_i\ra
\varphi_j)_{m:i\ra j}).$$
We specify a $\calC$-link $\varphi$ by a $\calC$-word $(\varphi_i)_{i\in |\varphi|}$ with
a set of isomorphisms each of which is composable only 
with its own inverses in the set. We call those isomorphisms also as  
\emph{the involutions of} the $\calC$-link $\varphi$. 

We can define the groupoid of $\calC$-links, whose isomorphism
$\kappa:\varphi\ra \psi$ between $\calC$-links is a pair
\pair{|\kappa|,\setii{\kappa_i}{i\in |\varphi|}}of an isomorphic
functor $|\kappa|:|\varphi|\ra |\psi|$ and isomorphisms
$$ \kappa_i:\varphi_i\ra \psi_{|\kappa|(i)}\qquad i\in |\varphi|$$
which makes the following commutative for $m:i\ra j$ in $|\varphi|$:
\begin{center}
  \includegraphics{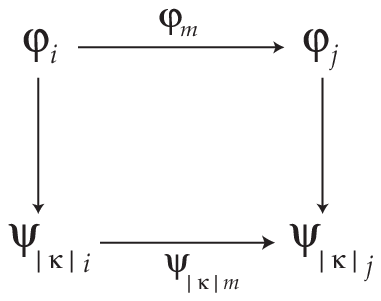}
\end{center}

$\calC$-links $\varphi_i$ has coprodct $\coprod_{i\in I}\varphi_i$ whose link
type is the coproduct of the link types of $\varphi_i$ and the functor
$ \coprod_{i\in I}\ra \calC$ is the direct sum of the functors $\varphi_i$.

There is a forgetful functor $U$ from the groupoid of $\calC$-links 
to that of $\calC$-words and isomorphisms defined
by throwing away the arrow parts:$  \varphi\mapsto U\varphi=(\varphi_i)_{i\in |\varphi|}$.
\emph{A link structure on} a $\calC$-word $w$ is 
a $\calC$-link $\varphi$ with $U\varphi=w$.



\section{Shells}
\subsection{Motivating example}
First we explain an example which motivate the definition of shells. 
The following is the process of blowing up the tetrahedron, first along edges and then
at the vertices. 

\begin{center}
 \includegraphics{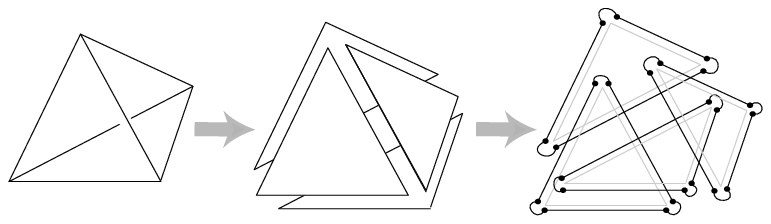}
\end{center}

The final result is can be described as in Fig. \ref{fig:tetrahedron-shell}.
\begin{figure}[htbp]
  \begin{center}
    \leavevmode
    \includegraphics{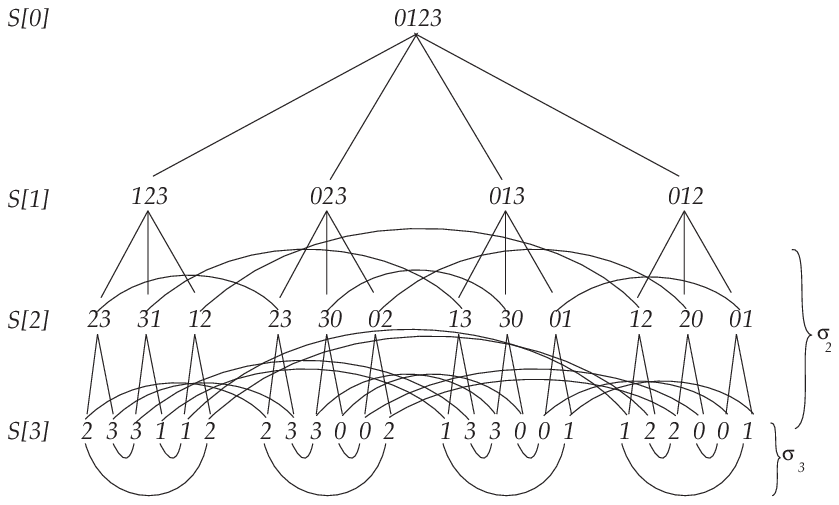}
    \caption{%
The tree representation of the $3$-shell representing the blow-up
structure of the tetrahedron. The two nodes are considered different
components even they have the same label. The vertex $0$ in the
tetrahedron appears as 6 $0$-components in $S[3]$ which are linked one another
by certain succession of isomorphisms  $\sigma^\bullet$.}
    \label{fig:tetrahedron-shell}
  \end{center}
\end{figure}

\subsection{Definition}
We formalize this combinatorial object as follows.

For nonnegative integers $n$, the notions of \emph{$n$-shells}, their isomorphisms
and closedness are defined by induction on $n$ as follows.

For $n=0$, an $n$-shell is simply a singleton set, regarded as the tree 
  consisting solely of the root. The unique map between singleton sets are the isomorphisms.
  Every $0$-shell is closed.

For $n=1$, an $n$-shell has the underlying tree $T_S$ of depth $\leq 1$.
  Its nodes of depth $1$ is regarded as $0$-shells. An isomorphism of $1$-shells is 
  simply a tree iso. Every $1$-shell is closed.

For $n\geq 2$, an $n$-shell consists of the following three data.
\begin{description}
\item[(Shell--1)] \emph{Its underlying tree $T_S$} of depth $n$, 
\item[(Shell--2)] For each node $x\in T_S[1]$, a closed $n-1$-shell $S^x$ with the
  underlying tree $T_S^x$. 
\end{description}

Note that not only the $n-1$-components but also every $i$-component $x$ with
  $i<n-1$ of an $n$-shell accompany a closed  $n-i$-shell
  $S^x$ and hence a link structure $\sigma^x$ on the word 
  $(S^y)_{y\in S^x[2]}$. From this it follows that for each $i>0$,
  there is an involution 
  $$\sigma_{n-i-2}:T[[i+2]]\ra T[[i+2]]$$
  and a partial involution 
  $$\sigma_{n-2}:T[[2]]\ra T[[2]]$$
  which commutes whenever composition is possible.

\begin{description}
\item[(Shell--3)] A link structure $\sigma_S$ on the word $(S^x)_{x\in T_S[2]}$ of 
  $n-2$-shells. 
\end{description}

An isomorphism $f:S \ra S'$ between $n$-shells  consists of the
following three data.
\begin{description}
\item[(Iso--1)] A tree isomorphism $|f|:T_S \ra T_{S'}$.  
\item[(Iso--2)] A shell isomorphism $f^x:S^x \ra S'^{|f|x}$ for each $x\in T_S[1]$.
\end{description}
Note that this induces $f^x:S^x \ra S'^{|f|x}$ for all $x\in T_S[[2]]$.
\begin{description}
\item[(Iso--3)] 
  A link isomorphism $\kappa:\sigma_S \ra \sigma_{S'}$ such that
  $$|\kappa|:=|f|_{T_S[2]}:|\sigma_S|=T_S[2]\ra |\sigma_{S'}|=T_{S'}[2]$$ 
  and, for each $x\in T_S[2]$, 
  $$\kappa_x:=f^x: \sigma_S(x)=S^x \ra \sigma_{S'}(|f|(x))=S'^{|f|(x)}.$$ 
\end{description}

An $n$-shell $S$ is \emph{closed} if the link structure $\sigma_S$ is closed,
which concludes the inductive defintion. 

An element of $T_S[i]$ of an $n$-shell $S$ is called \emph{of 
dimension $n-i$} and \emph{of codimension $i$}. We call an element of
dimension $k$ simply as \emph{a $k$-component}. 

A shell is called \emph{open} if it is not closed.

Two components are called \emph{linked} if the one is mapped to the other
by some involution $\sigma_i$. 

\subsection{Graphical representations}
\label{sec:graphical-rep}
\subsubsection{Tree representation}
A complete representation up to isomorphisms of the structure of an
$n$-shell is given by tree representation: We draw the underlying tree
and the tree iso's $\sigma^s$ by linking $x$ and $\sigma^s x$ for
$x$ in the domain of $\sigma^s$.

\subsubsection{Link representations}
Link representation allows us to 
focus our attention to crucial aspects of shells
by suitably neglecting inessential details of lower dimensional
components. For example, the closed 3-shell represented by the tree diagram for
the tetrahedron can be represented by the middle 
one in Figure \ref{fig:link-diagram} ignoring the $0$-components.
If necessary we can draw the information on the $0$-components as in the upper Figure \ref{fig:link-diagram}, which for 3-shells
has complete information. The lower Figure is another representation of link
diagram. For dimension grater than 4, it seems necessary to have
an appropriate method of concealing lower dimensional components. 

\begin{figure}
  \begin{center}
    \leavevmode
    \includegraphics{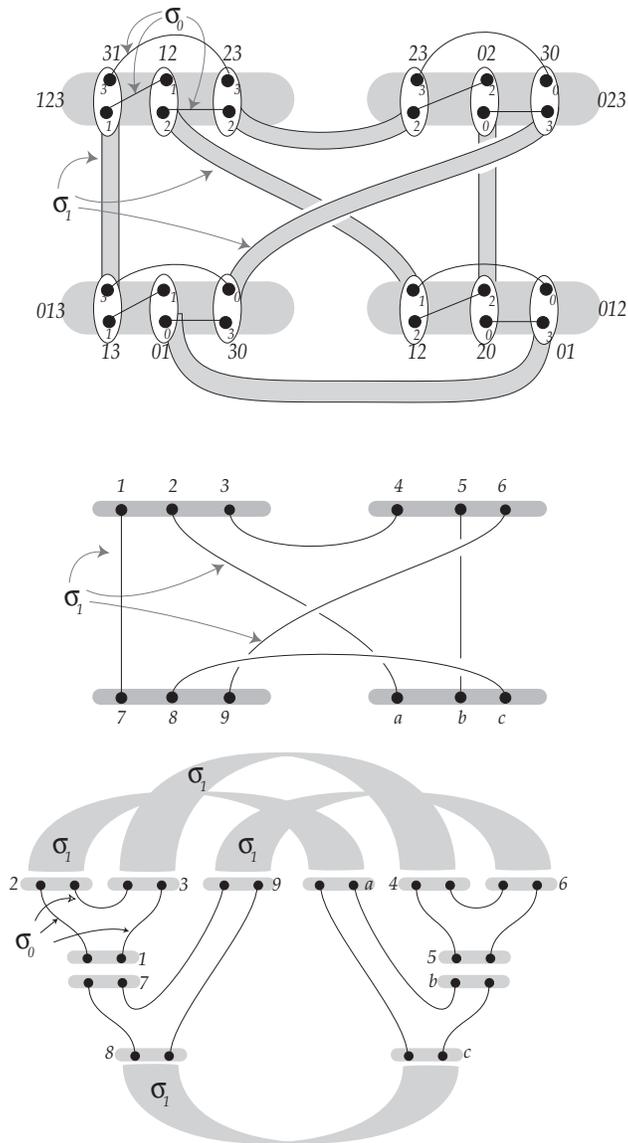}
    \caption{Link representations of the $3$-shell of the tetrahedron.
      The upper one keep all the
      information whereas the middle one, used frequently, neglects
      information on $0$-components. The lower one
      is another link diagram which give two dimensional forms to 2-components.}
    \label{fig:link-diagram}
  \end{center}
\end{figure}

Figure \ref{fig:C-3-shell.eps} draws an open $3$-shell in two ways.
\begin{figure}[htbp]
  \begin{center}
    \leavevmode
    \includegraphics{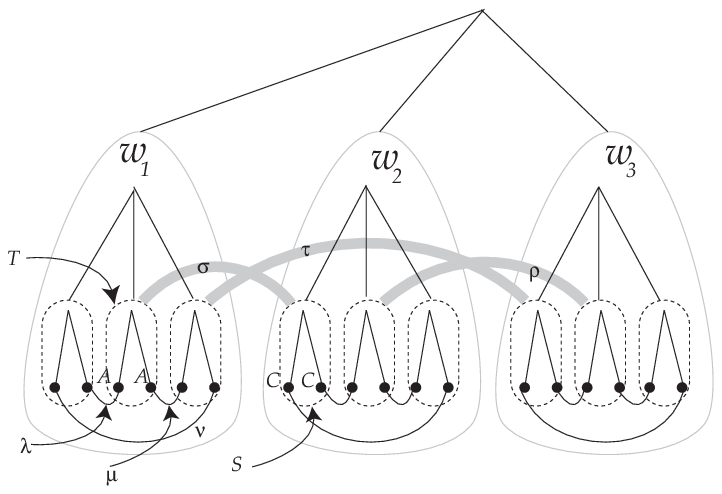}      \includegraphics{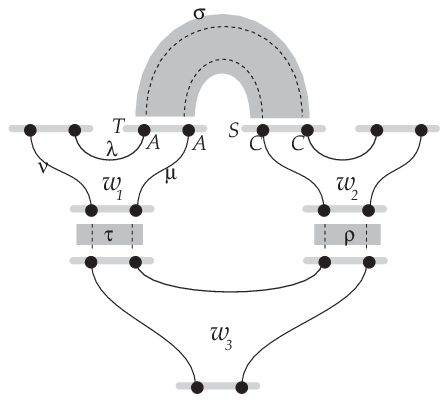}
    \caption{Descriptions of an open $3$-shell. In the left diagram, 
      $T_w[1]=\seti{1,2,3}$ and the $2$-shells $S^i$ are the $2$-shell of 
      the blow-up of a triangle. The involution $\sigma_3$ is represented by thin curves
      connecting the leaves and the involution $\sigma_2$ by thick curves.
      Note that, for example, the thick link $\sigma$ accompany an
      isomorphism between the $1$-shell $T$ and $S$. 
      The right diagram also describes the same pasting diagram, 
      which seems to be more appealing intuitively. 
}
    \label{fig:C-3-shell.eps}
  \end{center}
\end{figure}

\subsection{Closure of open shells}
\newcommand{\oT}{\overline{T}}
\newcommand{\oS}{\overline{S}}
From each open $n$-shell $S$ we construct a closed $n$-shell $\oS$, 
called \emph{its closure}. 

First define a new tree $\oT$ by adding new nodes to $T:=T_S$ as follows.
First we add a node $c_S$ of depth $1$ so that
$$ \oT[1]:=T[1]\;\coprod\;\seti{c_S} .$$
Let $T_{ext}[2]$ denotes the external indices of the link structure $\sigma_S$
and denote by $ T_{ext}[k]$ the set of nodes of depth $k$ which have antecedents 
in $T_{ext}[2]$. Define for $k\geq2$
$$ \oT[k]:=T[k]\;\coprod\;T_{ext}[k].$$ 
Note that $t\in T_{ext}[k]$ appears twice in the right hand side, so we use
$t$ for its occurrence in the first summand and $t^{*}$ for the second.
We say the nodes $c_s$ and $t^*$'s are \emph{new} whereas the nodes
$t$'s \emph{old}.

Now we define a shell structure on $\oT$. 

First we need to determine the $n-1$-shell $S^{c_S}$.  The only
missing information is the closed link structure on the word of
$n-3$-shells indexed by the set of external nodes
$T_{ext}[[3]]\subseteq T[3]$.  Let $\tau$ be the graph of [S--d] with
the vertex set $T[3]$. The external nodes in $T[3]$ are nothing but
the nodes of degree one and hence from each external node there is a
unique maximal chain with another external node as a terminal.  Hence
we obtained an involution $\mu$ on $T_{ext}[3]$. A chain connecting
$t$ and $\mu t$ are accompanied by tree iso's in the obvious way.

The Fig. \ref{fig:C-3-shell-boundary.eps} illustrates this construction for the open  
$3$-shell in Fig. \ref{fig:C-3-shell.eps} 

\begin{figure}[htbp]
  \begin{center}
    \leavevmode
    \includegraphics{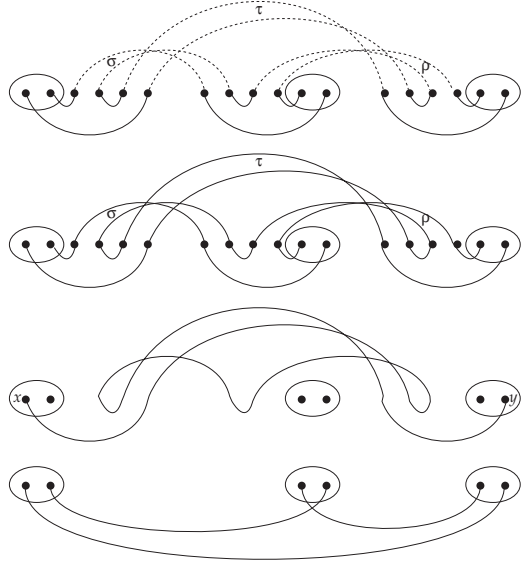}
    \caption{The construction of a link on the external nodes of depth $3$
      for the open $3$-shell in Fig.\ref{fig:C-3-shell.eps}.
      The six nodes in circles are the external ones.  The first
      figure describes two link types on $T_w[3]$, the dotted one from
      $\lambda$ and the solid one from $\kappa$.  The second one describes
      the graph $\Gamma$. The third one exhibits the path joining $x,y\in
L$.  The bottom one is the link structure for $S^{c_S}$.}
    \label{fig:C-3-shell-boundary.eps}
  \end{center}
\end{figure}

\newcommand{\overr}{{\overline{r}}}
For the root $\overr:=o_{\oT}$, the link $\sigma^{\overr}$ is defined
by extending the involutions of $\sigma^{o_T}$, and defining for the
external $t\in T[[2]]$ as the tree isomorphism $T^t \simeq T^{t^*}$ which maps
$u$ to $u^*$ for external $u\in T^t$. 

For a new node $t^*$ other than $c_S$, 
we define  $\sigma^{t^*}$ as the composition of the tree iso's:
$$ \oT^{u^*} \simeq T^u \stackrel{\sigma^t}{\simeq} T^v \simeq \oT^{v^*}$$
for $u,v\in T^t[[2]]$. 

It is obvious that $\oS$ is a closed shell. This completes the
construction of the closure $\oS$ of the shell $S$.

\section{Substantiation of shells}

\subsection{Labeling of shells}
Let $\Sigma$ be a set endowed with an involution $x \leftrightarrow
x^*$ called \emph{the conjugation} and a conjugation invariant grading
$\Sigma\ra \integer_+:=\seti{0,1,2,3,\cdots}$. 
The set of labels of grade $i$ is denoted by $\Sigma_i$.

We call such a set with a conjugation \emph{a labeling set}. 
We fix a labeling set $\Sigma$ in this section.
 
Let $S$ be an $n$-shell. 
A partial map $\lambda:T_S\ra \Sigma$ with the domain $|\lambda|\subseteq T_S$
is called \emph{a labeling of the shell $S$ with the label set $\Sigma$}  or
simply \emph{a labeling on $S$} if 
\begin{itemize}
\item it is graded in the sense that 
  the grade of $\lambda(x)$ is $k$ if $x$ is $k$-dimensional,
\item the domain $|\lambda|$ is descendant closed in the sense that
  if $\lambda$ is defined on $x$ then it is defined on the descendants of $x$,
\item it is compatible with the link structures, namely, two components
  have conjugate labels whenever they are linked.
\end{itemize}
\emph{The conjugate labeling  $\lambda^*$
of} a labeling $\lambda$ is defined by
$$ \lambda^*(s):=\lambda(s)^*.$$

A pair $(S,\lambda)$ of an $n$-shell and a labeling $\lambda$ is called 
\emph{a partial cell}. A partial cell $(S,\lambda)$ \emph{extends a partial cell }
$(S,\mu)$ if $|\mu|\subset |\lambda|$ and $\mu$ is the restriction of $\lambda$.

A labeling on an $n$-shell whose domain is the set of all the nodes of
depth $\geq n-k$ is called a $k$ dimensional labeling on $S$ over the
label set $\Sigma$, or simply \emph{a $k$-labeling on $S$ over $\Sigma$}.  If the
domain is the whole tree $T_S$, the pair $(S,\lambda)$ is called
\emph{an $n$-cell over $\Sigma$}.

Note that a partial labeling $\lambda$ restricts to labeling on
subsets of its domain.  In particular, if a component $s$ satisfies
$S^s\subseteq |\lambda|$, then $\lambda$ defines a total labeling on $S^s$. 
When $s$ is $i$-dimensional, the cell $(S^s, \lambda|S^s)$ is called
\emph{an $i$-face of} the partial cell $(S,\lambda)$.
For example, a $k$-labeling $\lambda$ defines boundary $i$-cells $(S^s,\lambda|S^s)$ 
for $i$-components $s$ with $i\leq k$.

\subsection{Cells and pasting diagrams}
We give the following names to pairs $(S,\lambda)$ of $n$-shells $S$
and $k$-labeling $\lambda$.

\vspace{2mm}
\begin{tabular}{c|c|c}
\hline
$S$ & $k=n$ & $k=n-1$ \\
\hline
closed & $n$-cell & $n$-frame \\
\hline
open & -- & $(n-1)$-pasting diagram\\
\hline
\end{tabular}
\vspace{2mm}

An $n$-cell $(S,\lambda)$ defines the $n$-frame
$(S,\lambda|S{[[1]]})$, called \emph{its boundary frame} or simply
\emph{the boundary}. The $n$-cell $(S,\lambda)$ is said to \emph{extend
the $n$-frame $(S,\lambda|S_{[[1]]})$}. 

We note also that every $n$-cell
$(S,\lambda)$ defines an $n$-pasting diagram $(\widetilde{S},\lambda)$
where $\widetilde{S}$ is an $n+1$-shell defined by
$\widetilde{S}[i]=S[i-1]$ for $i\geq0$, which we call \emph{the $n$-pasting
diagram with the only one $n$-component $(S,\lambda)$}.

\subsection{Examples}
The following shows examples
of 1-pasting diagram, 2-frame and 2-cell
$$    \includegraphics{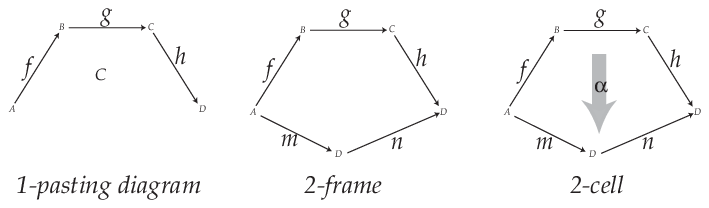}$$
in the usual way of describing cells of $2$-categories.

The Fig. \ref{fig:C-2-shell.eps} describes a $1$-pasting diagram $w$
using the representation method explained in \S \ref{sec:graphical-rep}.

\begin{figure}[htbp]
  \begin{center}
    \leavevmode
  \includegraphics{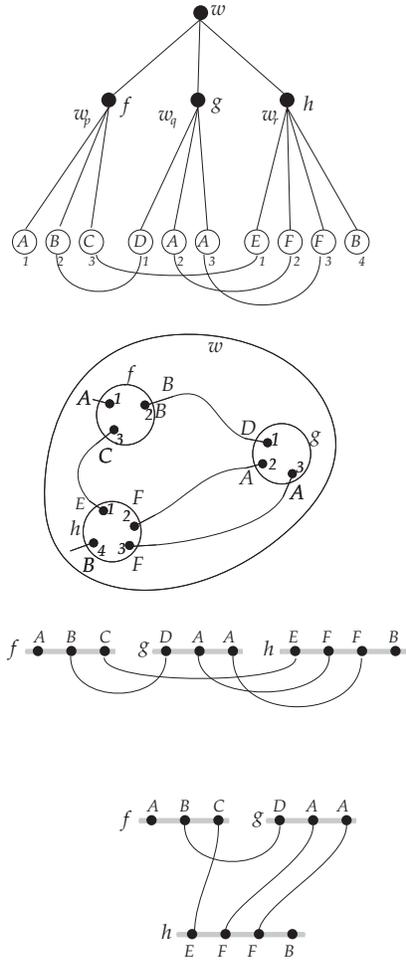}    
    \caption{A $1$-pasting diagram. 
      The symbols $f,g,h$ stand for $1$-hyperoperators,
      $A,B,C,D,E,F$ for $0$-hyperoperators with $A^*=F$, $B^*=D$, $C^*=E$,
      and $w,w_p,w_q,w_r,1,2,3,\cdots$ for the component of the underlying
      shell.
      }
\label{fig:C-2-shell.eps}
  \end{center}
\end{figure}

\subsection{The formal composites of pasting diagrams}
An $n$-pasting diagram describes an arbitrary combinatorially possible way 
of composing $n$-cells. The following theorem shows that an $n$-pasting diagram
uniquely defines a closed $n$-shell, called \emph{its formal composite}, 
which will be the underlying shell of the actual composite. 

\begin{thm}
Let $(S,\lambda)$ be an $n$-pasting diagram and $\oS$ be the closure of the
open $n$-shell $S$. Then
the labeling $\lambda$ extends uniquely to a labeling on $\oS{[[2]]}\bigcup 
S{[[1]]}$. 
In particular, it induces an $n-1$ labeling $\olam$ on the formal
composite $c(S)$ of $S$.
The $n$-frame $(c(S),\olam|c(S))$ is called \emph{the formal composite of }
the pasting diagram $(S,\lambda)$ and is denoted by $c(S,\lambda)$. 
\end{thm}

\begin{eproof}
Since every new component in $\oS$ are linked to some component in $S$, 
there is at most one extension of the labeling.
The existence of the extension is obvious for components of dimension less than
$n-3$. The compatibility of the labeling of new $n-3$-components
follows from the definition of $\osigma_{n-3}$. 
\end{eproof}

An $n+1$-cell $(\oS,\tilde\lambda)$ is called \emph{a composer of} the pasting diagram
$(S,\lambda)$ if $S|\tilde\lambda=\lambda$. The $n$-cell $(c(S),\tilde\lambda|c(S))$
is called \emph{the composite of the pasting diagram $(S,\lambda)$ by the $n+1$ cell 
$(\oS,\tilde\lambda)$} or simply \emph{a composite}. 

The lower is the formal composite of the upper $1$-pasting diagram, drawn in two ways. 
$$    \includegraphics{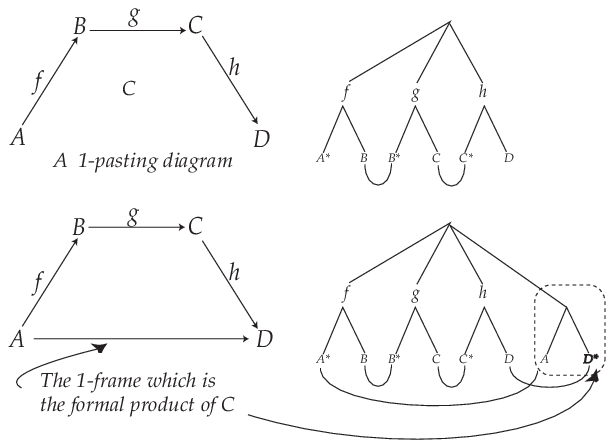}.$$

\section{Hypergraphs}
\subsection{Definition}
For $n\geq0$ we define the notions of $n$-hypergraphs,
the boundary operator and the coherence of labeling of $n+1$-shells over $\Sigma$.
by induction on $n$.
The set of coherent $n$-frames of an $n-1$-hypergraph
$\Sigma=\coprod_{0\leq i\leq n-1}\Sigma_i$ is denoted by
$\myframe_n{\Sigma}$. Element of $\Sigma_i$ is called \emph{a hyperoperator of
dimension $i$} or simply \emph{an $i$-hyperoperator.}

\begin{enumerate}
\item For $n=0$, every labeling set with $\Sigma=\Sigma_0$ is a $0$-hypergraph.
 Every labeling on a $1$-shell is defined to be coherent. 

\item A $1$-hypergraph is a labeling set $\Sigma=\Sigma_0\coprod
\Sigma_1$ with the boundary operator $\delta:\Sigma_1\ra
\myframe_1(\Sigma_0)$ which commutes with the conjugation
operators, i.e., $\delta(c^*)=\delta(c)^*$.  A labeling $\lambda$ on a $2$-shell $S$
over $\Sigma$ is called coherent if for every
$t\in T_S[1]\bigcap |\lambda|$ , the $1$-frame $(T_S^t,\lambda|T_S^t[1])$ coincides with
the boundary $1$-frame $\delta(\lambda(t))$ of the 
$1$-hyperoperator $\lambda(t)\in \Sigma_1$.

\item For $n\geq 2$, an $n$-hypergraph is a labeling set
$$ \Sigma=\coprod_{0\leq i\leq n}\Sigma_i$$
with boundary operators
$$ \delta_i:\Sigma_i\ra \myframe_i(\Sigma)$$
for $1\leq i\leq n$ such that 
$$(\Sigma\setminus \Sigma_n,\delta_1,\cdots,\delta_{n-1}) $$
is an $n-1$-hypergraph. A labeling $\lambda$ on an $n+1$-shell $S$
is called coherent if for every $i$-component $t$ of $S$, the
$i$-frame $(S^t,\lambda|S^t[[1]])$ coincides with the boundary of the
$i$-hyperoperator $\lambda(t)$. 
\end{enumerate}

We denote an $n$-hypergraph as 
$$(\coprod_{0\leq i\leq n}\Sigma_i,\delta_1,\cdots,\delta_{n}) $$
or simply as $(\Sigma,\delta:\Sigma\ra\myframe(\Sigma))$. 
An $\omega$-hypergraph is a labeling set $\Sigma$
with the boundary operator $ \Sigma_m\ra \myframe_m(\Sigma)$ for all
natural number $m$.

\emph{Hereafter all labelings over hypergraphs will be assumed to be coherent. }

\subsection{Category of hypergraphs}
Let $\calH,\calH'$ be $n$-hypergraphs. A graded map $f:\calH\ra \calH'$ commuting
with the conjugation is called \emph{a hypergraph map } if
it commutes with the boundary map, in the sense that for every $i$-edge $t$ of $\calH$,
$$ \delta(f(t))= f_i\delta(t),$$
where $f_i:\myframe_i\calH\ra \myframe_i\calH'$ is defined by composing
$f$ to the labeling map. We denote by $\hypergraph_n$ the category of
$n$-hypergraphs and hypergraph maps.

\section{The monad structure}
Let $(\calH_0,\cdots,\calH_n,\delta_1,\cdots,\delta_n)$ 
be an $n$-hypergraph. 
Denote by $\PD_n\calH$ the set of $n$-pasting diagrams over $\calH$.
Since the formal composite of a coherent pasting diagram 
is a coherent frame, we have 
$$\comp:\PD_n\calH\ra\myframe_n\calH=\myframe_n(\PD\calH),$$ 
where $\PD$ is  a  labeling set 
with $(\PD\calH)_i=\calH_i$ for $i<n$ and $(\PD\calH)_n=\PD_n\calH$,  
with the conjugation on $\PD_n\calH$ given by the conjugation of the labeling.

\subsection{Multiplication on $\PD$}
This induces an endofunctor $\PD$ of $\hypergraph_n$ defined by
\begin{eqnarray*}
\PD\calH&=&(\calH_0,\cdots,\calH_{n-1},
\PD_n\calH,\\
&&\delta_1,\cdots,\delta_{n-1},\comp:\PD_n\calH\ra\myframe_n(\PD\calH)).
\end{eqnarray*}
This endofunctor has a monad structure whose multiplication
$$  \mu:\PD\PD\calH \ra \PD\calH$$
is defined as follows. Since  $(\PD\PD\calH)_i=\calH_i$ for $i<n$, we define $\mu_i$ to be
the identity map. So we need only to define
$$ \mu_n:(\PD\PD\calH)_n=\PD_n(\PD\calH)\ra \PD_n\calH.$$

Let $(S,\lambda)$ be an $n$-pasting diagram over $\PD\calH$, where $S$ is an open $n+1$-shell.
For $t\in T_S[1]$, the labeling $\lambda(t)\in \PD_n\calH$ is an $n$-pasting diagram $(U_t,\lambda_t)$
over $\calH$ whose boundary is the $n$-frame$(S^t,\lambda|S^t[[1]])$, i.e., $S^t$ is the external
components of $T_{U_t}$. 

Define now an $n+1$-shell $V$. Its underlying tree is defined by 
$$ T_V[i]=\coprod_{t\in T_S[1]}T_{U_t}[i]$$
for $i\geq 1$.  For each $\in T_{U_t}[i]$ for $i>0$, we define 
$i$-shell $V^t:=U^s_t$. The link structure for the root of $T_V$ is defined by joining those for the
roots of $T_{U_t}$'s and that for the root of $T_S$ transferred to the link structure for the
words of shells indexed by the external nodes of $T_{U_t}$'s. 
Fig. \ref{fig:monad} illustrates this construction.
\begin{figure}[htbp]
  \begin{center}
    \leavevmode
    \includegraphics{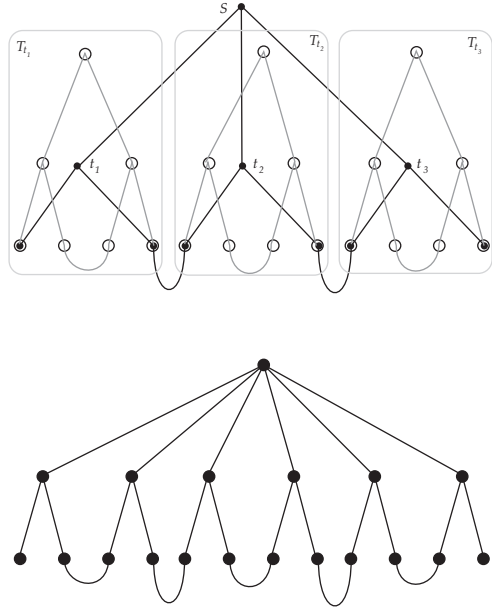}
    \caption{The lower one is the composition of the upper pasting diagram of the pasting diagram. }
    \label{fig:monad}
  \end{center}
\end{figure}

The unit natural transformation $\eta:\calH\ra \PD\calH$ is given 
by defining $\eta_n(s)$ to be the $n$-pasting diagram consisting solely of $s$.

It is straightforward to see the following.

\begin{theorem}
The triple $(\PD,\mu,\eta)$ is a monad on the category of hypergraphs.
\end{theorem}

\subsection{Submonads of $(\PD,\mu,\alpha)$}
By restricting the shapes of the pasting diagrams, we obtain
a few submonads of $\PD$ which are equally important.

An $n$-pasting diagram $\varphi=(S,\lambda)$ over $\calH$ defines a graph 
a graph $\Gamma_\varphi$ with the vertex set $T_S[2]$. An element \seti{i,j} 
is its edge if and only if there is an $n$-component $s\in T_S[1]$ whose
link structure $\sigma^s$ on $T^s_S[2]\subseteq T_S[3]$ connects a child of $i$ with
that of $j$. For example, the $3$-cell of Fig. \ref{fig:C-3-shell.eps} is as follows:  
\begin{center}
  \includegraphics{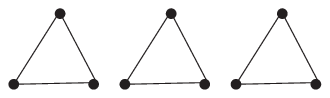}
\end{center}

A pasting diagram $\varphi$ is called \emph{acircuit} if the graph $\Gamma_\varphi$
has no circuits. It is called \emph{connected} if the graph $\Gamma\varphi$ is connected.
Denote by $\psdmtree\calH$ the collection of connected acircuit pasting diagrams of $\calH$.

\newcommand{\sign}{\mbox{sgn}}
Suppose now that every $n-1$-hyperoperator $x$ of an $n$-hypergraph
$\calH$ has signature $\sign(x)\in \seti{-1,1}$ with
$\sign(x^*)=-\sign(x)$. Let $\varphi=(S,\lambda)$ be an $n$-pasting diagram.
We give the direction to the graph $\Gamma_\varphi$ 
by $i\ra j$ if there are linked $p\in \child(i)$ and $q\in \child(j)$ with $\lambda{}p$ positive and $\lambda{}q$ negative.  We say the pasting diagram $\varphi$ is
\emph{acyclic} if this graph has no cycles.  Let $\psdmposet\calA$ be
the collection of connected acircuit pasting diagrams of $\calX$.

We have
\begin{prop}
Both $\psdmtree$ and $\psdmposet$ are submonad of $\PD$.
\end{prop}
Thus we obtain two monads $(\psdmtree,\mu,\eta)$ and
$(\psdmposet,\mu,\eta)$.

\newcommand{\hgraph}{{\cal H}}
\newcommand{\hgcat}{{\cal HG}raph}

\section{Strict Hypercategories} 
\subsection{Definition}
Let $n\geq 1$.
\emph{A strict $n$-hypercategory} is a $\PD$-algebra 
$(\calH,\alpha:\PD\calH\ra \calH)$ where $\calH$ is an $n$-hypergraph.

Each submonad of $\PD$, defines a variant of 
\emph{strict $n$-hypercategory}. We call a $\psdmtree$-algebra 
\emph{an acircuit hypercategory}. Similarly we define \emph{
an acyclic hypercategory} 
when the hyperoperators have parities.

The monad multiplication $\mu$ defines \emph{the free
$n$-hypercategory $\PD\calH$ generated by an $n$-hypergraph $\calH$}.

We first give examples of $1$-hypercategories and then
illustrate the usefulness of $2$-hypercategories by 
explicitly representing the rewriting process of  the
Lafont's interaction combinator as a $2$-pasting diagram of
a strict $2$-hypergraph.

\subsection{Examples of $1$-hypercategories}
\subsubsection{Classical simple logic}
Let $\calH_0$ be a set of propositional variables together with
their negations. The conjugation is the negation. 
Fix a truth assignment of the variables. 
Write $\proves P_1,\cdots,P_m$ when
at least one of $P_1,\cdots,P_m$ is true. 
If we define $\calH_1$ the set of finite sets \seti{P_1,\cdots,P_m}
of propositional variables with $\proves P_1,\cdots,P_m$, then 
we have a strict acircuit $1$-hypercategory with the multiplication 
given by the following cut rule:
$$ 
\Frac{\proves P,P_1,\cdots,P_n\quad\proves P^*,Q_1,\cdots,Q_m}{
\proves P_1,\cdots,P_n,Q_1,\cdots,Q_m}.
$$

\subsubsection{Categories as acircuit hypercategories}
Let $\calC$ be a category. 
We define simply 
$$ A^*=\overline{A}\qquad \overline{A}^*=A.$$

An arrow $f:A\ra B$ defines an arrow $f$ with 
$$ \delta f=(A,\overline{B}).$$
An acircuit pasting diagram is simply a linear diagram of composable sequence of arrows.
\begin{prop}
  There is a bijection between categories with objects $\calC_0$ 
and acircuit hypercategories over $\calC_0\coprod \overline{\calC_0}$ whose
hyperedges have boundaries of the form $(A,\overline{B})$.  
\end{prop}

\subsubsection{Multicategories}
Let $\calM$ be a multicategory in the sense of Lambek. 
For each object $A$ we prepare its conjugate $\overline{A}$ 
and define $\calH_0$ to be the collection of objects and their conjugates.

An arrow $ \Gamma:A_1,\cdots,A_n\ra B$ corresponds to an arrow
$\varphi_\Gamma$ with the boundary $(A_1,\cdots,A_n,\overline{C})$.
Then by the associativity of the cut rule, we obtain a strict acircuit
$1$-hypercategory. In fact it can be seen easily that there is a bijection
between multicategories with objects $\calH_0$ and acircuit
hypercategories over $(\calC_0\coprod \overline{\calC_0})$.

\subsection{Lafont's interaction net}
\label{sec:2d}

Let $\calH_0$ be the singleton set $\seti{a}$. 
Let $\calH_1$ consist of $1$-hyperoperators
$\seti{0,\epsilon,s,+,\times,\delta}$, whose boundaries are
as follows: 
$$
\begin{array}{c|c|c|c|c|c}
\hline
0 & \epsilon & s & + & \times & \delta \\
\hline
a & a & aa & aaa & aaa & aaa\\
\hline
\end{array}.
$$
We denote these as
\begin{center}
  \includegraphics{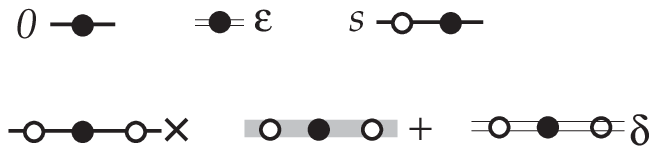}
\end{center}

Let $\calH_2$ be the set of 2-hyperoperators 
whose boundaries are described in Fig. \ref{fig:lafont-2-hyperoperator}.

\begin{figure}[htbp]
  \begin{center}
    \leavevmode
  \includegraphics{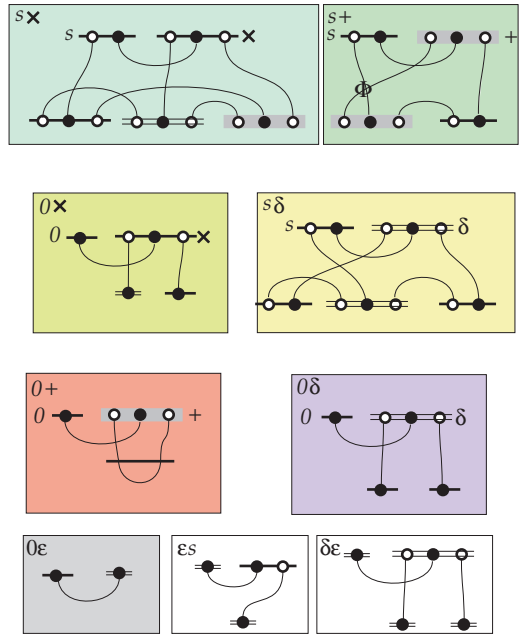}    
    \caption{The set of 2-hyperoperators}
    \label{fig:lafont-2-hyperoperator}
  \end{center}
\end{figure}

For example, the $s+$ 2-cell is described by tree form as:
\begin{center}
  \includegraphics{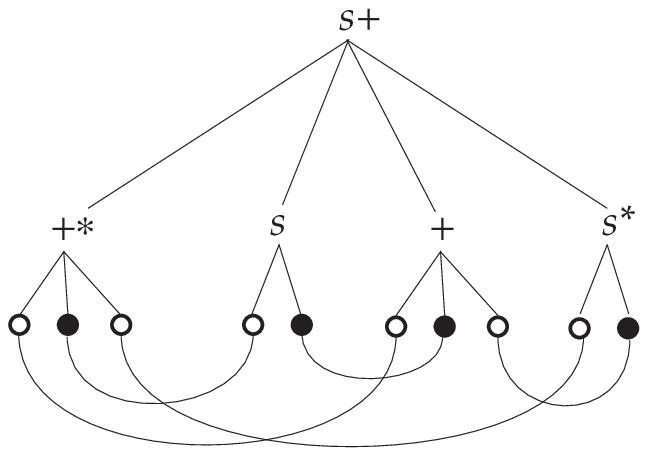}
\end{center}

The following is an example of pasting diagram of the lafont's $2$-hypergraph.
\begin{center}
  \includegraphics{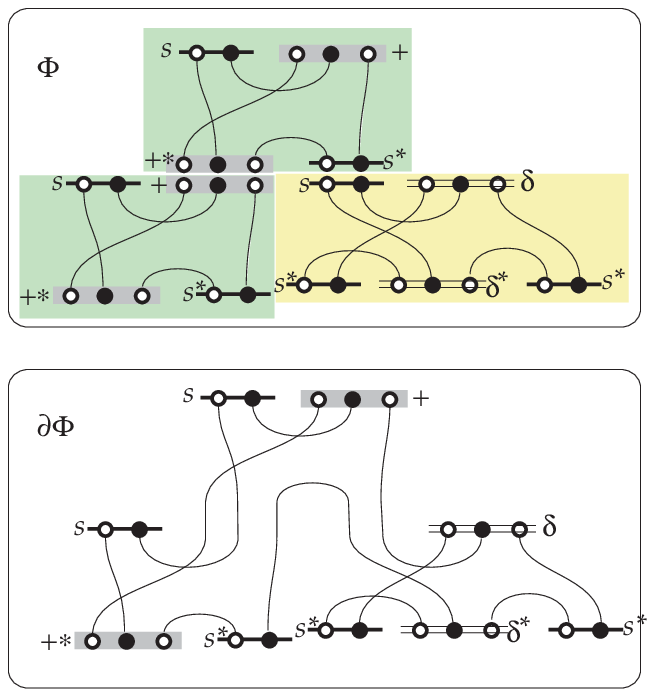}
\end{center}

The Figure \ref{fig:lafont-addition.eps} is the pasting
diagram which corresponds to the sequence of actions of interaction net
which calculates $2 \times 2=4$.
\begin{figure}[htbp]
  \begin{center}
    \leavevmode
\includegraphics{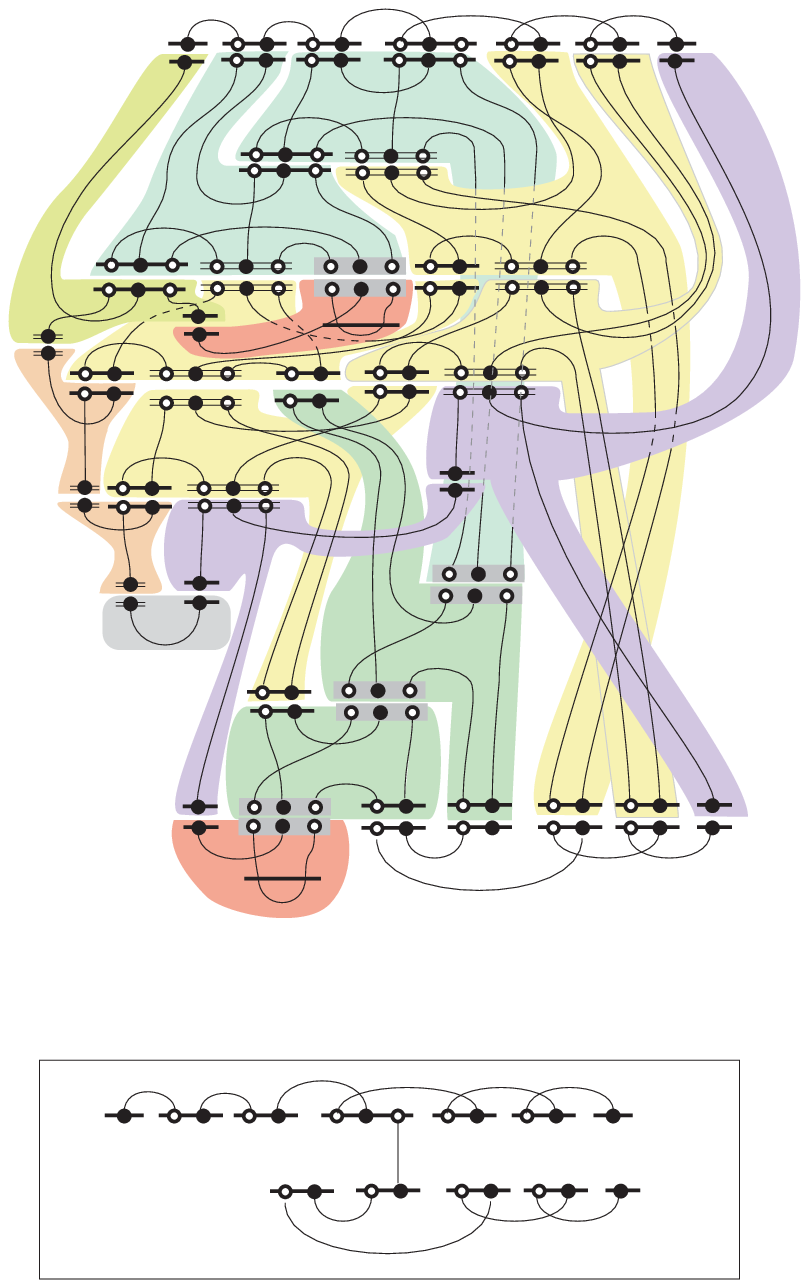}    
    \caption{The pasting diagram of the interactions calculating $2\times 2=4$. }
    \label{fig:lafont-addition.eps}
  \end{center}
\end{figure}

\section{Hypercategory}
We sketch a formulation of weak $\omega$-hypercategory based on hypergraphs.

\subsection{Definition}
We assume that the underlying $\omega$-hypergraph $\calH$ satisfy
the following conditions. (i) Every hyperoperator $f$ has the parity
$\sign(f)\in \seti{1,-1}$ with $\sign(f^*)=-\sign(f)$.
A partial $m+1$-cell $(S,\lambda)$ over $\calH$ is called 
\emph{pure} if the labels of $m-1$-components are of the same sign.
(ii) The boundary of hyperoperators are pure.

\emph{A hypercategory} is an $\omega$-hypergraph
$\calH=((\calH_i)_{i=0,1,\cdots},\delta)$ with certain elements
called \emph{universal} are singled out and satisfy the following
conditions.
\noindent\textbf{(H${}_1$)} Every pure $n$-pasting diagram $C$ has a universal
composer $U$. The composite of $C$ with respect to the composer $U$ is called \emph{
universally composed $C$}. Similarly every pure $n+1$-frame $F$ has a universal
$n+1$-cell $U$ with $\partial U=F$.
\\
\noindent\textbf{(H${}_2$)} Every universal $m$-cell $U$ has a universal 
cell $U^{\dagger}$ of the same sign with the frame conjugate to $\partial U$, called
\emph{a transpose}\footnote{This should not be confused with $U^*$ which has opposite parity
and exists for all $U$. Moreover there are usually more than one transposes}.\\
\noindent\textbf{(H${}_3$)} If a pure $n$-pasting diagram is universal in the sense that 
the labeling of $n$-components are universal, then its composites universally composed are universal. 

The condition (\textbf{H${}_2$}) 
replaces the role of that involving \emph{balancedness} in
Baez-Dolan's definition. In fact we can easily show the following.

\begin{prop} Let $\calH$ be a hypercategory and and suppose 
an $n$-pasting diagram $(S,\lambda) $ has two composites $C,C'$.
\noindent\textrm{(i)} If $C$ is  universally composed, then there is an 
$n+1$-cell $M$ whose frame has the underlying $n+1$-shell, denoted by
$F_{n+1}(C^*,C')$, whose components of dimension $\leq n$ belongs either to $C$ or to $C'$
and the involution on $n-1$-components is given by the identity of the $n$-shell $S$. \\
\noindent\textrm{(ii)} If $C'$ and $C$ are both universal, then 
then there is a universal $n+1$-cell $C$ whose boundary is $F_{n+1}(C^*,C')$.
\end{prop}
\begin{eproof}
  Just take the following $n+1$-pasting diagram. Its $n+1$-components
are a transpose of the universal composer of $(S,\lambda)$ composing
$C$ and the $n+1$-cell giving the composite $C'$. The $n-1$ involution
the identity on $S{[[1]]}$. Then its composite gives the asserted $n+1$ cell.
\end{eproof}

A hypercategory $\calH$ is called \emph{of dimension $n$} or simply an
$n$-hypercategory if every cell of dimension greater than $n$ is
universal.  It is called \emph{$m$-weak} if every pure pasting
diagram of dimension greater than $m$ has a unique composite.

\subsection{Hypercategories over an $n$-hypergraph}
Baez-Dolan give a method of restricting class of $n$-categories by restricting
the type of shells. Some of their procedure can be described by using 
prototype $n$-hypergraph. 

Let $\Sigma$ be an $n$-hypergraph. An $\omega$-hypergraph $\calH$ with
a hypergraph map $\varphi:\calH_{[n]}\ra \Sigma$ is called \emph{an
$\omega$-hypergraph over $\Sigma$}. A hypercategory with the
underlying hypergraph over $\Sigma$ is said to \emph{be of type
$\Sigma$.}

\subsection{$0$-weak $1$-hypercategory}
To show some aspects of arguments in our formulation, we show that usual category
is obtained from $0$-weak $1$-hypercategory. 
\label{section:type}
Let $\Sigma_0=\seti{a,a^*},\Sigma_1=\seti{b,b^*}$ and $\delta
b=(a^*,a)$.  A $0$-weak $1$-hypercategory over $\Sigma$ corresponds to
a usual category in the following way.

First of all, pure 1-pasting diagrams are nothing but composable
sequences of arrows and by the 0-weakness, they have the unique
composite, called \emph{the composition} which is universally
composed.

\begin{lem}
  Let $f,g:A\ra B$. If there is a universal 2-cell $u:f\ra g$, then $f=g$. 
\end{lem}
\begin{eproof}
Note that $u$ is a universal composer of $f$ regarded as a 1-pasting diagram. 
Since $u$ is universal, there is another universal $u^\dagger:g\ra f$ which together
gives a universal composite $f\ra f$, which is another universal composer of the
$f$ regarded as a $1$-pasting diagram and hence must coincide with $u$. 
In particular $f=g$.
\end{eproof}
\begin{prop}
  The composition is associative. 
\end{prop}
\begin{eproof}
For simplicity, let us prove $f\circ (g\circ h)=(f\circ g)\circ h$. The pasting diagram
$$\includegraphics{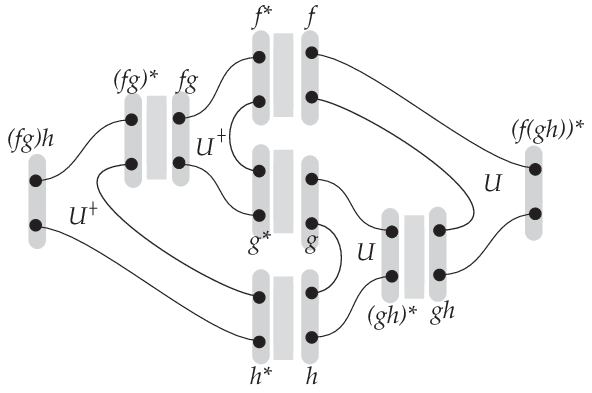}$$
gives a universal 2-cell $f\circ (g\circ h) \ra (f\circ g)\circ h$, whence the associativity follows from the above lemma.
\end{eproof}
Each object $A\in \calH_0$, regarded as a $0$-pasting diagram has universal
composer which we call \emph{quasi-identities} temporally. 

\begin{prop}
\textrm{(i)} Let $u:A\ra A$ be a quasi-identity. If $f:A\ra B$, then $f\circ u=f$.  Similarly
if $g:B\ra A$, then $u\circ g=g$.  
(ii) For each $A$, there is only one quasi-identity.
\end{prop}
\begin{eproof}
Let $0_u$ be the universal composer of the 1-pasting diagram $u$. Then the 
2-pasting diagram 
$$\includegraphics{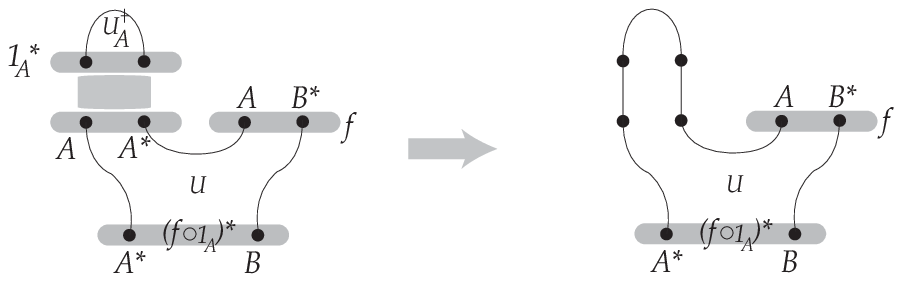}$$
gives a composite $f\circ u \ra f$ universally composed, whence by the lemma $f\circ u=f$.
If $u,v$ are quasi-identities, then $u=u\circ v=v$.
\end{eproof}
Similarly multicategory can be identified with a class of $0$-weak $1$-hypercategories.

\end{document}